\documentclass[12pt]{amsart}
\usepackage{amsmath, amssymb, latexsym, amsthm, mathrsfs, color, mathtools}
\usepackage[hyphens]{url}

\usepackage[top=2in, bottom=1.5in, left=1in, right=1in]{geometry}

\newcommand{\nc}{\newcommand}

\nc{\thusfar}{\par\bigskip\centerline{\my{--- Edited thus far ---}}\par\bigskip}
\nc{\lei}{\le^\oo}
\nc{\card}[1]{\left|#1\right|}
\nc{\medcard}[1]{\biggl|\,#1\,\biggr|}
\nc{\smallcard}[1]{|\,#1\,|}
\nc{\bds}{bidirectional $\roth$-scale}
\nc{\bbN}{\mathbb{N}}
\nc{\beq}{\begin{eqnarray*}}\nc{\eeq}{\end{eqnarray*}}
\nc{\mbq}{\mb{?}}
\nc{\mb}[1]{{\mbox{\textbf{#1}}}}
\nc{\nop}{$\times$}
\nc{\fbn}{\!\!\fbox{\!\nop\!}\!\!}
\nc{\yup}{\checkmark}
\nc{\forces}{\Vdash}
\nc{\name}[1]{\dot{#1}}
\nc{\tf}{\my{FINISHED THUS FAR}}
\nc{\FU}{Fr\'echet--Urysohn}
\nc{\gs}{$\gamma$~space}
\nc{\Ga}{\Gamma}\nc{\Om}{\Omega}\nc{\Lam}{\Lambda}
\nc{\smallbinom}[2]{\begin{psmallmatrix} #1\\ #2 \end{psmallmatrix}}
\nc{\bgamma}{\smallbinom{\Om}{\Ga}}
\newcommand{\two}{\{0,1\}}
\nc{\productive}[2]{\bigl(#1,\allowbreak #2\bigr)^\x}
\nc{\Sel}{\mathsf{S}}
\nc{\sset}[2]{\{\,#1 : #2\,\}}
\nc{\smb}[1]{{\!\!\mb{#1}\!\!}}
\nc{\medset}[2]{{\biggl\{\,#1 : #2\,\biggr\}}}
\nc{\smallmedset}[2]{{\bigl\{\,#1 : #2\,\bigr\}}}
\nc{\set}[2]{{\left\{\,#1 : #2\,\right\}}}
\nc{\seq}[2]{{\la\, #1 : #2\,\ra}}
\nc{\eseq}[1]{#1_1, \allowbreak #1_2, \allowbreak\dotsc} 
\nc{\cube}{(\Cantor)^\bbN}
\nc{\Match}{\op{Match}}
\nc{\concat}[1]{\hat{\phantom{a}}\langle #1\rangle}
\nc{\poset}{\mathbb{P}}
\nc{\fn}[1]{{\op{Fn}(#1\times\w,2)}}
\nc{\linadd}{\op{linadd}}
\nc{\nonprod}{\non^\x}
\nc{\alephes}{{\aleph_0}}
\nc{\my}[1]{\marginpar{\textcolor{red}{***}}\textcolor{red}{#1}}
\nc{\later}[1]{{\color{green} #1}}
\nc{\BTs}[1]{{\color{green} #1 (BT)}}
\nc{\Cp}{\op{C}_\mathrm{p}}
\nc{\Bp}{\op{B}_p}
\nc{\Pa}[8]{\bibitem{#1} {#2}, \emph{#3}, {#4} \textbf{#5} ({#6}), {#7}--{#8}.}
\nc{\tPa}[5]{\bibitem{#1} {#2}, \emph{#3}, {#4}, to appear.}
\nc{\sPa}[4]{\bibitem{#1} {#2}, \emph{#3}, {#4}, submitted.}
\nc{\Bc}[9]{\bibitem{#1} {#2}, \emph{#3}, in: \textbf{#4} (#5), #6 #7, #8--#9.}
\nc{\fD}{\mathfrak{D}}
\nc{\fX}{\mathfrak{X}}
\nc{\Onbd}{\Op_{\mathrm{nbd}}} 
\nc{\Omnb}{\Om_{\mathrm{nbd}}} 
\nc{\od}{\mathfrak{od}}
\nc{\Setting}[7]{\xymatrix@R=4pt@C=7pt{#1\ar@{-}[r]&#2\ar@{-}[r]&#3\\&#4\ar@{-}[u]\\
#5\ar@{-}[uu]\ar@{-}[r] & #6\ar@{-}[u]\ar@{-}[r] & #7\ar@{-}[uu]}}
\nc{\mx}[1]{\begin{matrix}#1\end{matrix}}
\nc{\plim}{p\txt{-}\lim}
\nc{\Bgp}{{\Z^\bbN}}
\nc{\Cgp}{{{\Z_2}^\bbN}}
\nc{\Cite}[1]{\textbf{[#1]}}
\nc{\Next}[1]{{#1^+}}
\nc{\cFin}{\mathrm{cF}}
\nc{\intvl}[2]{{[#1(#2),\allowbreak #1(#2\!+\!1))}}
\nc{\Bdd}{\mathbf{B}}
\nc{\Dfin}{\mathfrak{D}_\mathrm{fin}}
\nc{\grbl}{{\mbox{\textit{\tiny gp}}}}
\nc{\bbP}{\mathbb{P}}
\nc{\BOfat}{\B_{\Om_{\mathrm{fat}}}}
\nc{\Bgood}{\B_{\mathrm{good}}}
\nc{\compactN}{\cl{\mathbb{N}}}
\nc{\blocks}[2]{\op{cl}_{#2}(#1)}
\nc{\blocksplus}[2]{\op{cl}^+_{#2}(#1)}
\nc{\arx}[1]{\texttt{http://arxiv.org/math/#1}}
\nc{\bq}{\begin{quote}}
\nc{\eq}{\end{quote}}
\nc{\cl}[1]{\overline{#1}}
\nc{\CH}{the Continuum Hypothesis}
\nc{\MA}{Martin's Axiom}
\nc{\Bfat}{\B_\mathrm{fat}}
\nc{\inv}{^{-1}}
\nc{\Cantor}{{\two^\bbN}}
\nc{\bP}{\mathbf{P}}
\nc{\bof}{\op{\fb}}
\nc{\dof}{\op{\fd}}
\nc{\bofF}{\bof(\cF)}
\nc{\sr}[3]{\underset{\mbox{#3}}{\mbox{#1}}}
\nc{\gp}{\binom{\Om}{\Ga}}
\nc{\gpsmall}{\mbox{$\gp$}}
\nc{\gig}{\gimel}
\nc{\gns}{\sone(\Om,\gig)}
\nc{\nsr}[2]{#1}
\nc{\Srg}{{\mathbb{S}}}
\nc{\Srgs}{{\mathbb{S}^*}}
\nc{\NN}{{\bbN^{\bbN}}}
\nc{\ZN}{{\Z^{\bbN}}}
\nc{\NNup}{{\bbN^{\uparrow\bbN}}}
\nc{\Pof}{\op{P}}
\nc{\PN}{{\Pof(\bbN)}}
\nc{\roth}{{[\bbN]^{\mbox{\tiny $\infty$}}}} 
\nc{\Fin}{[\bbN]^{\text{$<\!\!\infty$}}} 
\nc{\FinSeqs}{\bbN^{\text{$<\!\!\infty$}}} 
\nc{\ici}{[\bbN]^{ \infty, \infty}}
\nc{\Inc}{{\compactN^{\uparrow\bbN}}}
\nc{\powInc}[1]{{\big(\Inc\big)^{#1}}}
\nc{\powFin}[1]{{\big(\Fin\big)^{#1}}}
\nc{\powPN}[1]{{\big(\PN\big)^{#1}}}
\nc{\NcompactN}{{\compactN^\bbN}}
\nc{\Uarrow}{\smash{\big\uparrow}}
\nc{\LE}{\preccurlyeq}
\nc{\GE}{\succcurlyeq}
\nc{\op}{\operatorname}
\nc{\im}{\op{im}}
\nc{\Span}{\op{span}}
\nc{\maxfin}{\op{maxfin}}
\nc{\ran}{\op{range}}
\nc{\iso}{\cong}
\nc{\Madd}{{\M}^\star}
\nc{\cI}{\mathcal{I}}
\nc{\cJ}{\mathcal{J}}
\nc{\scrA}{\mathscr{A}}
\nc{\scrB}{\mathscr{B}}
\nc{\scrC}{\mathscr{C}}
\nc{\scrD}{\mathscr{D}}
\nc{\scrF}{\mathscr{F}}
\nc{\scrK}{\mathscr{K}}
\nc{\A}{\forall}
\nc{\B}{\mathrm{B}}
\nc{\cB}{\mathcal{B}}
\nc{\bB}{\mathbf{B}}
\nc{\BS}{\mathbf{B}(\mathcal{S})}
\nc{\BF}{\mathbf{B}(\mathcal{F})}
\nc{\BU}{\mathbf{B}(\mathcal{U})}
\nc{\cSp}{\mathcal{S}^+}
\nc{\cFp}{\mathcal{F}^+}
\nc{\cUp}{\mathcal{U}^+}

\nc{\BG}{\B_\Ga}
\nc{\BL}{\B_\Lambda}
\nc{\BT}{\B_\Tau}
\nc{\BTstar}{\B_{\Tau^*}}
\nc{\BO}{\B_\Om}
\nc{\DO}{\cD_\Om}
\nc{\KO}{\cK_\Om}
\nc{\CG}{C_\Ga}
\nc{\CL}{C_\Lambda}
\nc{\CT}{C_\Tau}
\nc{\CTstar}{C_{\Tau^*}}
\nc{\CO}{C_\Om}
\nc{\COgp}{C_{\Om^{\grbl}}}
\nc{\CLgp}{C_{\Lambda^{\grbl}}}
\nc{\BOgp}{\B_{\Om}^{\grbl}}
\nc{\BLgp}{\B_{\Lambda^{\grbl}}}
\nc{\sfC}{\mathsf{C}}
\nc{\sfD}{\mathsf{D}}
\nc{\bD}{\mathbf{D}}
\nc{\Tau}{\mathrm{T}}
\nc{\cA}{\mathcal{A}}
\nc{\cK}{\mathcal{K}}
\nc{\cD}{\mathcal{D}}
\nc{\cF}{\mathcal{F}}
\nc{\cS}{\mathcal{S}}
\nc{\cT}{\mathcal{T}}
\nc{\cG}{\mathcal{G}}
\nc{\cY}{\mathcal{Y}}
\nc{\J}{\mathcal{J}}
\nc{\cL}{\mathcal{L}}
\nc{\cM}{\mathcal{M}}
\nc{\cN}{\mathcal{N}}
\nc{\cH}{\mathcal{H}}
\nc{\cO}{\mathcal{O}}
\nc{\Op}{\mathrm{O}}
\nc{\rmA}{\mathrm{A}}
\nc{\rmF}{\mathrm{F}}
\nc{\rmB}{\mathrm{B}}
\nc{\rmD}{\mathrm{D}}
\nc{\rmP}{\mathrm{P}}
\nc{\cC}{\mathcal{C}}
\nc{\cP}{\mathcal{P}}
\nc{\bbQ}{\mathbb{Q}}
\nc{\bbR}{\mathbb{R}}
\nc{\cU}{\mathcal{U}}
\nc{\Un}{\bigcup}
\nc{\cV}{\mathcal{V}}
\nc{\cW}{\mathcal{W}}
\nc{\Z}{{\mathbb Z}}
\nc{\Impl}{\Rightarrow}
\long\def\forget#1\forgotten{\marginpar{\textcolor{green}{Forgetting...}}}
\nc{\ft}{\mathfrak{t}}
\nc{\fb}{\mathfrak{b}}
\nc{\fc}{\mathfrak{c}}
\nc{\fd}{\mathfrak{d}}
\nc{\fg}{\mathfrak{g}}
\nc{\oo}{\infty}
\nc{\fr}{\mathfrak{r}}
\nc{\fk}{\mathfrak{k}}
\nc{\bidi}{\mathfrak{bidi}}
\nc{\fu}{\mathfrak{u}}
\nc{\fh}{\mathfrak{h}}
\nc{\fp}{\mathfrak{p}}
\nc{\fj}{\mathfrak{j}}
\nc{\fs}{\mathfrak{s}}
\nc{\w}{\omega}
\nc{\x}{\times}
\nc{\Iff}{\Leftrightarrow}

\nc{\nin}{\notin}
\nc{\cat}{\hat{\ }}
\nc{\sub}{\subseteq}
\nc{\spst}{\supseteq}
\nc{\sm}{\setminus}
\nc{\as}{\subseteq^*}
\nc{\les}{\le^*}
\nc{\leinf}{\le^{\infty}}
\nc{\leS}{\le_S}
\nc{\leF}{\le_{\mathcal{F}}}
\nc{\leU}{\le_{\mathcal{U}}}
\nc{\rest}{\restriction}

\nc{\la}{\langle}
\nc{\ra}{\rangle}
\nc{\E}{\exists}
\nc{\dom}{\op{dom}}
\nc{\cov}{\op{cov}}
\nc{\add}{\op{add}}
\nc{\cof}{\op{cof}}
\nc{\cf}{\op{cf}}
\nc{\non}{\op{non}}
\nc{\unif}{\op{non}}
\nc{\COV}{\op{COV}}
\nc{\ADD}{\op{ADD}}
\nc{\COF}{\op{COF}}
\nc{\NON}{\op{NON}}
\nc{\impl}{\to}
\nc{\Lp}{\mathcal{L_\p}}
\nc{\Wlog}{without loss of generality}
\newtheorem{thm}{Theorem}
\nc{\bthm}{\begin{thm}} \nc{\ethm}{\end{thm}}
\newtheorem{prop}[thm]{Proposition}
\nc{\bprp}{\begin{prop}} \nc{\eprp}{\end{prop}}
\newtheorem{fact}[thm]{Fact}
\nc{\bfct}{\begin{fact}} \nc{\efct}{\end{fact}}
\newtheorem{prob}[thm]{Problem}
\nc{\bprb}{\begin{prob}} \nc{\eprb}{\end{prob}}
\newtheorem{lem}[thm]{Lemma}
\nc{\blem}{\begin{lem}} \nc{\elem}{\end{lem}}
\newtheorem{claim}[thm]{Claim}
\nc{\bclm}{\begin{claim}} \nc{\eclm}{\end{claim}}
\newtheorem{cor}[thm]{Corollary}
\nc{\bcor}{\begin{cor}} \nc{\ecor}{\end{cor}}
\newtheorem{conj}[thm]{Conjecture}
\nc{\bcnj}{\begin{conj}} \nc{\ecnj}{\end{conj}}
\theoremstyle{definition}
\newtheorem{defn}[thm]{Definition}
\nc{\bdfn}{\begin{defn}} \nc{\edfn}{\end{defn}}
\newtheorem{obs}[thm]{Observation}
\nc{\bobs}{\begin{obs}} \nc{\eobs}{\end{obs}}
\theoremstyle{remark}
\newtheorem{rem}[thm]{Remark}
\nc{\brem}{\begin{rem}} \nc{\erem}{\end{rem}}
\newtheorem{cnv}[thm]{Convention}
\nc{\bcnv}{\begin{cnv}} \nc{\ecnv}{\end{cnv}}
\newtheorem{exam}[thm]{Example}
\nc{\bexm}{\begin{exam}} \nc{\eexm}{\end{exam}}
\nc{\bpf}{\begin{proof}} \nc{\epf}{\end{proof}}
\nc{\be}{\begin{enumerate}}
\nc{\ee}{\end{enumerate}}
\nc{\bi}{\begin{itemize}}
\nc{\bimy}{\my{\begin{itemize}}
\nc{\eimy}{\end{itemize}}}
\nc{\itm}{\item}
\nc{\ei}{\end{itemize}}
\nc{\Subsection}[1]{\goodbreak\subsection*{#1}}

\nc{\sone}{\mathsf{S}_1}
\nc{\sfin}{\mathsf{S}_\mathrm{fin}}
\nc{\ufin}{\mathsf{U}_\mathrm{fin}}
\nc{\Split}{\mathsf{Split}}

\nc{\gone}{\mathsf{G}_1}    \nc{\gfin}{\mathsf{G}_\mathrm{fin}}

\nc{\ed}{
	
\subsection*{Acknowledgments}
We thank Micha\l{} Machura and Jialiang He for helpful discussions that helped us
understand Pawlikowski's proof better.
The first named author thanks the second for his hospitality during 2016,
where the present proof of Hurewicz's Theorem was obtained.
The second named author thanks the first for
his hospitality beyond the conference \emph{Frontiers of Selection Principles} (Warsaw, 2017),
where the main breakthroughs leading to the present proof of Pawlikowski's Theorem were made.
A preliminary version of this paper was quoted in a survey of Aurichi and Dias~\cite{AD}, and this led
to a ``flood'' of requests by colleagues to see these notes, and this urged us to complete this paper.
We thank our colleagues for their enthusiasm.

\end{document}
}

\title[The Menger and Rothberger games]{Conceptual proofs of the Menger and Rothberger games}

\author[P. Szewczak]{Piotr Szewczak}
\address{Piotr Szewczak,
Institute of Mathematics, Faculty of Mathematics and Natural Science College of Sciences, Cardinal Stefan Wyszy\'nski University in Warsaw, Warsaw, Poland,
and
Department of Mathematics, Bar-Ilan University, Ramat Gan, Israel
}
\email{p.szewczak@wp.pl}
\urladdr{http://piotrszewczak.pl}

\author[B. Tsaban]{Boaz Tsaban}
\address{Boaz Tsaban,
Department of Mathematics, Bar-Ilan University, Ramat Gan, Israel}
\email{tsaban@math.biu.ac.il}
\urladdr{http://math.biu.ac.il/~tsaban}

\begin{document}

\begin{abstract}
	We provide conceptual proofs of the two most fundamental theorems
	concerning topological games and open covers:
	Hurewicz's Theorem concerning the Menger game,
	and Pawlikowski's Theorem concerning the Rothberger game.
\end{abstract}

\maketitle



\section{Introduction}

Topological games form a major tool in the study of topological properties and their relations to
Ramsey theory, forcing, function spaces, and other related topics.
Excellent surveys~\cite{SchOPiT, AD}, and a comprehensive bibliography for further inspection~\cite{Tel},
are available.
At the heart of the theory of \emph{selection principles}~\cite{wiki},
covering properties are defined by the ability to diagonalize, in canonical ways,
sequences of open covers. Each of these covering properties has an associated two-player game.
Often the nonexistence of a winning strategy for the first player in the associated game is equivalent
to the original property, and this forms a strong tool for establishing results concerning the
original property.

Allowing few potential exceptions, all proofs of results of this type use, as black boxes,
two fundamental theorems.
The original proofs of the fundamental theorems are technical and difficult to digest, and the
lack of their understanding has undoubtedly veiled deep insights
that necessitate variations in the original proofs.
We present here new, conceptual proofs of these theorems.
These proofs build on the earlier proofs, but are intuitive and can be grasped
without the need of technical verifications.

\section{The Menger game}

A topological space has \emph{Menger's property} $\sfin(\Op,\Op)$ if, for each sequence
of open covers, $\cU_1,\cU_2,\dots$, we can select finite sets $\cF_1\sub\cU_1$, $\cF_2\sub\cF_2$,\dots whose union $\Un_n\cF_n$ covers the space.
The symbol $\Op$ in this notation indicates that we are provided with open covers, and need to obtain an open cover.
These properties were also considered for additional classes of open covers, in the realm of 
selection principles~\cite{wiki}.

\emph{Menger's game} $\gfin(\Op,\Op)$
is a game for two players, Alice and Bob, with an
inning per each natural number $n$.
In each inning, Alice picks an open cover of the ambient space and Bob selects finitely many
members from this cover.
Bob wins if the sets he selected throughout the game cover
the space. If this is not the case, Alice wins.

If Alice does not have a winning strategy in the game $\gfin(\Op,\Op)$, then
$\sfin(\Op,\Op)$ holds. The converse implication is a deep theorem of Hurewicz~\cite[Theorem~10]{Hure25}.
This theorem is often used within selection principles, and has applications
in diverse contexts, such as
D-spaces~\cite{AurD} and additive Ramsey theory~\cite{AlgSelRT}.
We present here a conceptual proof of this theorem.
The proof uses the same initial simplifications as in Scheepers's proof of Hurewicz's Theorem~\cite[Theorem~13]{coc1}.
We add one further simplification that makes calculations easier, and
an appropriate notion that goes through induction, and thus
eliminates the necessity to track the history of the game.

\bthm[Hurewicz]
\label{thm:hur}
Let $X$ be a space satisfying $\sfin(\Op,\Op)$ space. Then Alice does not have a winning strategy in the
game $\gfin(\Op,\Op)$.
\ethm
\bpf
Fix an arbitrary strategy for Alice in the game $\gfin(\Op,\Op)$.

If there is a play in which Bob covers the space after finitely many steps, then we are done.
Thus, we assume that, in no position, a finite selection suffices, together with the earlier selections, to cover the space.
Since the space $X$ satisfies $\sfin(\Op,\Op)$, it is Lindel\"of.
By restricting Bob's moves
to countable subcovers of Alice's covers, we may assume that
Alice's covers are countable.
Given that, we may assume that Alice's covers are \emph{increasing}, that is, of the form
$\{U_1,U_2,\dotsc\}$ with $U_1\sub U_2\sub\dotsb$, and that Bob selects a \emph{single} set in
each move. Indeed, given a countable cover $\{U_1,U_2,\dotsc\}$, we can
restrict Bob's selections to the form $\{U_1,U_2,\dotsc,U_n\}$, for $n\in\bbN$.
Since Bob's goal is just to cover the space, we may pretend that Bob is provided covers
of the form
\[
\{U_1,U_1\cup U_2,U_1\cup U_2\cup U_3,\dotsc\},
\]
and if Bob selects an element $U_1\cup\dotsb\cup U_n$, he replies to Alice with the legal
move $\{U_1,U_2,\dotsc,\allowbreak U_n\}$. If Bob manages to cover the space, this is not due to the unions.

Finally, we may assume that for each reply $\{U_1,U_2,\dots\}$ (with $U_1\sub U_2\sub\dotsb$)
of Alice's strategy to a move
$U$, we have $U=U_1$.
Indeed, we can transform the given cover into the cover
$\{U, U\cup U_1, U\cup U_2,\dots\}$. If Bob chooses
$U$, we provide Alice with the answer $U_1$, and if he chooses
$U\cup U_n$, we provide Alice with the answer $U_n$.
Since Bob has already chosen the set $U$,
its addition in the new strategy does not help covering more points.

With these simplifications, Alice's strategy is identified with a tree of open sets,
as follows:
Alice's initial move is an open cover $\{U_{1},U_{2},\dots\}$.
If Bob replies $U_{n}$, then Alice's move is $\{U_{n,1},U_{n,2},\dots\}$.
In general, if Bob replies $U_\sigma$, for $\sigma\in\bbN^k$, then Alice's move
is an increasing open cover
\[\cU_\sigma:=\{U_{\sigma(1),\dotsc,\sigma(k),1},U_{\sigma(1),\dotsc,\sigma(k),2},\dots\},\]
with $U_\sigma=U_{\sigma(1),\dotsc,\sigma(k),1}$.

The proof will reduce to the following concept.

\bdfn
A countable cover $\cU$ of a space $X$ is a \emph{tail cover} if the set of
intersections of cofinite subsets of $\cU$ is an open cover of $X$.
\edfn

Equivalently, a cover $\{U_1,U_2,\dots\}$ is a tail cover if
the family
\[
\Bigl\{\,\bigcap_{n=1}^\infty U_n,\bigcap_{n=2}^\infty U_n,\dotsc\,\Bigr\}
\]
of intersections of cofinal segments of the cover is an open cover.

\blem
Let $n$ be a natural number. Define
$\cV_n:=\Un_{\sigma\in\bbN^n}\cU_\sigma$. Then the family $\cV_n$ is a tail cover of $X$.
\elem
\bpf
The proof is by induction on $n$.

The open cover $\cV_1=\cU_{()}$ is increasing, and thus the set of
cofinite intersections is again $\cV_1$,
an open cover of $X$.

Let $n$ be a natural number. For brevity, enumerate $\cV_n=\{V_1,V_2,\dotsc\}$, and
\[
\cV_{n+1}=\Un_{k=1}^\infty\{V^k_1,V^k_2,\dotsc\},
\]
where
\[
V_k=V^k_1\sub V^k_2\sub\dotsb.
\]
We assume, inductively, that the family $\cV_n$ is a tail cover of $X$.
Let $\cV$ be a cofinite subset of $\cV_{n+1}$.
For each natural number $k$,
let $m_k$ be the minimal natural number with $V^k_{m_k}\in \cV$.
Then $\bigcap(\cV\cap \{V^k_1,V^k_2,\dotsc\})=V^k_{m_k}$
for all $k$, and $m_k=1$ for all but finitely many natural numbers $k$.
Let $I:=\set{k}{m_k=1}$, a cofinite subset of $\bbN$. We have
\[
\bigcap\cV=\bigcap_{k\in\bbN}(\cV\cap\{V^k_1,V^k_2,\dotsc\})=
\bigcap_{k\in\bbN} V^k_{m_k}=
\bigcap_{k\in I}V_k\cap \bigcap_{k\in \bbN\sm I}V^k_{m_k},
\]
Since $\cV_n$ is a tail cover, the set $\bigcap_{k\in I}V_k$ is open.
The remaining part is a finite intersection of open sets. Thus, the
set $\bigcap\cV$ is open.

Let $x\in X$. Since $\cV_n$ is a tail cover,
the set $I:=\sset{k}{x\in V_k}$ is cofinite.
For $k\in \bbN\sm I$, let $m_k$ be the minimal natural number with $x\in V^k_{m_k}$.
Then
$x\in \bigcap_{k\in I}V_k\cap \bigcap_{k\in \bbN\sm I}V^k_{m_k}$ and,
as we saw, the latter set is an intersection of a cofinite subset of
the family $\cV_{n+1}$.
\epf

For each $n$, let $\cV_n'$ be the set of intersections of cofinite subsets of $\cV_n$.
Applying the property $\sfin(\Op,\Op)$ to the sequence $\cV_1',\cV_2',\dotsc$,
Bob obtains cofinite sets $\cW_n\sub\cV_n$ such that $X=\Un_n\bigcap\cW_n$.
In the $n$-th inning, Alice provides Bob with a cover that is an infinite subset of the family
$\cV_n$. Since the family $\cW_n$ is cofinite in $\cV_n$,
Bob can choose an element $V_n\in\cV_n\cap\cW_n$.
Then $X=\Un_n V_n$, and Bob wins.

This completes the proof of Hurewicz's Theorem.
\epf

To treat Rothberger's game, we need a result slightly stronger than Hurewicz's.
The original proof of the following result, due to Pawlikowski~\cite[Lemma~1]{Paw94},
is much more technical than the combination of the present proofs of
Theorem~\ref{thm:hur} and Corollary~\ref{cor:paw1}.

We recall that Menger's property $\sfin(\Op,\Op)$ is preserved by countable unions:
Given a countable union of Menger spaces, and
a sequence of open covers, we can split the sequence of covers
into infinitely many disjoint subsequences,
and use each subsequence to cover one of the given Menger spaces.

\bcor[Pawlikowski]
\label{cor:paw1}
Let $X$ be a space satisfying $\sfin(\Op,\Op)$. For each strategy for Alice in the game
$\gfin(\Op,\Op)$, there is a play according to this strategy,
\[
(\cU_1,\cF_1,\cU_2,\cF_2,\dotsc),
\]
such that for each point $x\in X$ we have $x\in\Un\cF_n$ for infinitely many $n$.
\ecor
\bpf
We apply a reduction of Scheepers, originally used to prove the analogous theorem for
the game considered in Section~\ref{sec:roth}~\cite[Theorem~3]{OpPar}.

The product space $X\x\bbN$, a countable union of Menger spaces, satisfies $\sfin(\Op,\Op)$.
We define a strategy for Alice in the game $\gfin(\Op,\Op)$, played on the space $X\x\bbN$.
Let $\cU$ be Alice's first move in the original game.
Then, in the new game, her first move is
\[
\tilde\cU:=\set{U\x \{n\}}{U\in\cU, n\in\bbN}.
\]
If Bob selects a finite set $\tilde\cF\sub\tilde\cU$, we take the set
\[
\cF:=\set{U\in\cU}{\text{there is }n\text{ with }U\x\{n\}\in\tilde\cF}
\]
as a move in the original game.
Then Alice replies with a cover $\cV$, and we continue in the same manner.
By Hurewicz's Theorem, there is a play
\[
(\tilde\cU_1,\tilde\cF_1,\tilde\cU_2,\tilde\cF_2,\dotsc)
\]
in the new game, with $\Un_n\tilde\cF_n$ a cover of $X\x \bbN$.
Consider the corresponding play in the original strategy,
\[
(\cU_1,\cF_1,\cU_2,\cF_2,\dotsc).
\]
Let $x\in X$. There is a natural number $n_1$ with $(x,1)\in\Un\tilde\cF_{n_1}$.
Then $x\in\Un\cF_{n_1}$.
The set
\[
F:=\set{k\in\bbN}{\text{there is }U\text{ with } U\x\{k\}\in\Un_{i=1}^{n_1}\tilde\cF_i}
\]
is finite.
Let $m$ be a natural number greater than all elements of the set $F$.
There is a natural number $n_2$ with $(x,m)\in\Un\tilde\cF_{n_2}$.
Then $x\in\Un\cF_{n_2}$, and $n_1<n_2$. Continuing in a similar manner,
we see that $x\in\Un\cF_n$ for infinitely many $n$.
\epf

\brem
A cover of a space is \emph{large} if each point is covered by infinitely many
members of the cover.
Let $\Lam$ be the family of all large covers of an ambient space $X$.
We have $\sfin(\Op,\Op)=\sfin(\Lam,\Lam)$~(\cite[Corollary~5]{coc1}, \cite[Theorem~1.2]{coc2}).
With some initial simplifications of the considered strategies,
Corollary~\ref{cor:paw1} implies that a topological space $X$ satisfies
$\sfin(\Lam,\Lam)$ if and only if Alice does not have a winning strategy in the
corresponding game $\gfin(\Lam,\Lam)$. In fact, these results are essentially
identical.
\erem

\section{The Rothberger game}
\label{sec:roth}

The definitions of \emph{Rothberger's property} $\sone(\Op,\Op)$
and the corresponding game $\gone(\Op,\Op)$ are similar to those of
$\sfin(\Op,\allowbreak\Op)$
and
$\gfin(\Op,\Op)$,
respectively,
but here we select \emph{one} element from each cover.
Here too, if Alice does not have a winning strategy then the space satisfies
$\sone(\Op,\Op)$. The converse implication was established by Pawlikowski~\cite[p.~279]{Paw94}, improving considerably over partial results of
Galvin~\cite[Corollary~4]{Galvin77} and Rec\l{}aw~\cite[Corollary~2]{Rec94}.
We provide a conceptual proof of Pawlikowski's Theorem.
We first isolate an argument in Pawlikowski's proof, that does not involve
games.

\blem[Pawlikowski]
\label{lem:paw2}
Let $X$ be a space satisfying $\sone(\Op,\Op)$.
Let $\cF_1,\cF_2,\dotsc$ be nonempty finite families of open sets such that, for each point
$x\in X$, we have $x\in\Un\cF_n$ for infinitely many $n$.
Then there are elements $U_1\in\cF_1, U_2\in\cF_2,\dotsc$ such that
the family $\{U_1,U_2,\dotsc\}$ covers the space $X$.
\elem
\bpf
For each natural number $n$, let $\cU_n$ be the family of all intersections of $n$ open sets
taken from distinct members of the sequence $\cF_1,\cF_2,\dotsc$.
Then $\cU_n$ is an open cover of $X$.

By the property $\sone(\Op,\Op)$, there are elements $V_1\in\cU_1, V_2\in\cU_2,\dotsc$
that cover the space $X$. Extend the set $V_1$ to an element of some family $\cF_n$.
We can extend $V_2$ to an element of some \emph{other} family $\cF_n$, and so on.
We obtain a selection of at most element from each family $\cF_n$, that covers $X$.
We can extend our selection to have an element from each family $\cF_n$.
\epf

For a natural number $k$ and families of sets $\cU_1,\dotsc,\cU_k$, let
\[
\cU_1\wedge\dotsb\wedge \cU_k:=\set{U_1\cap\dotsb\cap U_k}{U_1\in\cU_1,\dotsc, U_k\in\cU_k}.
\]

\bthm[Pawlikowski]
\label{thm:paw3}
Let $X$ be a space satisfying $\sone(\Op,\Op)$.
Then Alice does not have a winning strategy in the game $\sone(\Op,\Op)$.
\ethm
\bpf
Fix an arbitrary strategy for Alice in the Rothberger game $\gone(\Op,\Op)$.
Since $\sone(\Op,\Op)$ spaces are Lindel\"of,
we may assume that each cover in the strategy is countable.
Let $\FinSeqs$ be the set of finite sequences of natural numbers.
We index the open covers in the strategy as
\[
\cU_\sigma=\{U_{\sigma,1},U_{\sigma,2},\dotsc\},
\]
for $\sigma\in\FinSeqs$, so that
$\cU=\{U_1,U_2,\dots\}$ is Alice's first move,
and for each finite sequence $k_1,\dotsc,k_n$ of natural numbers,
$\cU_{k_1,\dotsc,k_n}$ is Alice's reply to the position
\[
(\cU,U_{k_1},\cU_{k_1},U_{k_1,k_2},\cU_{k_1,k_2},\dotsc,U_{k_1,\dots,k_{n}}).
\]
For finite sequences $\tau,\sigma\in\bbN^n$, we write $\tau\le\sigma$ if $\tau(i)\le\sigma(i)$ for all $i=1,\dots,n$.
We define a strategy for Alice in the Menger game $\sfin(\Op,\Op)$.
Alice's first move is $\cU$, her first move in the original strategy. Assume that
Bob selects a finite subset $\cF$ of $\cU$.
Let $m_1$ be the minimal natural number with $\cF\sub\{U_{1},\dotsc,U_{m_1}\}$.
Then, in the Menger game, Alice's response is the joint refinement
$\cU_1\wedge\dotsb\wedge \cU_{m_1}$.
Assume that Bob chooses a finite subset $\cF$ of this refinement.
Let $m_2$ be the minimal natural number such that $\cF$ refines all sets $\{U_{i,1},\dots,U_{i,m_2}\}$,
for $i=1,\dotsc,m_1$. Then Alice's reply is the joint refinement
$\bigwedge_{\tau\le (m_1,m_2)} \cU_\tau$.
In general, Alice provides a cover of the form
$\bigwedge_{\tau\le \sigma} \cU_\tau$, for $\sigma\in\FinSeqs$, Bob selects a finite family
refining all families $\{U_{\tau,1},\dots,U_{\tau,m}\}$ for $\tau\le \sigma$, whit the minimal natural number $m$,
and Alice replies  $\bigwedge_{\tau\le (\sigma,m)} \cU_\tau$.

By Pawlikowski's Theorem (Corollary~\ref{cor:paw1}), there is a play
\[
(\cU,\cF_1,\bigwedge_{k_1\le m_1}\cU_k,\cF_2,\bigwedge_{(k_1,k_2)\le (m_1,m_2)}\cU_{k_1,k_2},\dotsc),
\]
according to the new strategy, such that every point of the space is covered
infinitely often in the sequence $\Un\cF_1,\Un\cF_2,\dotsc$.
By Lemma~\ref{lem:paw2}, we can pick one element from each set $\cF_n$ and cover
the space.
There is $k_1\le m_1$ such that the first picked element is a subset of $U_{k_1}$.
There is $k_2\le m_2$ such that the second picked element is a subset of $U_{k_1,k_2}$,
and so on.
Then the play
\[
(\cU,U_{k_1},\cU_{k_1},U_{k_2,k_2},\dots)
\]
is in accordance with Alice's strategy in the Rothberger game, and is won by Bob.
\epf

\ed

\appendix

\section{An alternative proof of Pawlikowski's theorem on Menger's game}

\my{Probably, to be removed in the final version}

We provide below an alternative proof for Pawlikowski's Theorem~\ref{???}.
This proof may be applicable to types of covers that cannot be treated by the
above proof.

\bthm[Pawlikowski]
Let $X$ be an $\sfin(\Op,\Op)$ space. Then Alice does not have a winning strategy
in the game $\gfin(\Lam,\Lam)$.
\ethm
\bpf
We first provide a slightly more technical proof of Hurewicz's theorem, and then show that it
provides the present theorem.
Given a strategy for Alice in the game $\gfin(\Lam,\Lam)$, we first restrict
the moves of Bob by removing, from each cover in the strategy, the finitely many
sets chosen on the way to this cover. The covers remain large. This way, it suffices
that, in his play, Bob covers each point of the space $X$ infinitely often, and he need not
worry that this may happen by finitely many open sets.
We may assume, further, that the covers have no finite subcover. Indeed, if the covers with
finite subcovers are dense in the strategy tree, Bob easily wins. And if not, we can restrict Bob to reach a node in the tree where, henceforth, there are no finite subcovers.
The first few moves contribute finitely many sets, and largeness is not affected by that.

Again, we may assume that the given covers are countably infinite and increasing (and have
no finite subcover).
Let $\cU_\sigma$, for $\sigma\in\FinSeqs$, be the open covers in Alice's strategy.
For increasing covers $\cU=\{U_1,U_2,\dotsc\}$ and $\cV=\{V_1,V_2,\dotsc\}$, define
\[
\cU\wedge\cV:=\{U_1\cap V_1,U_2\cap V_2,\dotsc\}.
\]
Then the cover $\cU\wedge\cV$ is increasing, too.

For $\tau,\sigma\in\bbN^n$, we write $\tau\le\sigma$ if $\tau(i)\le\sigma(i)$ for all $i=1,\dots,n$.
For each sequence $\sigma\in\FinSeqs$, define
\[
\cV_\sigma=\bigwedge_{\tau\le\sigma}\cU_\tau,
\]
and enumerate this cover, in its natural order,
as $\{V_{\sigma\cat 1},V_{\sigma\cat 2},\dotsc\}$.
Thus,
\[
V_{\sigma\cat n}=\bigcap_{\tau\le\sigma}U_{\tau\cat n}
\]
for all $n$.
The tree $\set{V_\sigma}{\sigma\in\FinSeqs}$ defines a new strategy for Alice,
whose covers are finer than the original ones. It thus suffices to prove that
the new strategy is not winning.

In summary, we may assume that
the original covers are countable, increasing, and getting finer as $\sigma$ gets larger,
and we have the following observation.

\bclm
\label{clm:mon}
For all sequence $\tau\le\sigma$, and natural numbers $m\le n$, we have
$U_{\sigma\cat m}\sub U_{\tau\cat n}$.
\eclm
\bpf
Since the open covers are increasing, we have
\[
U_{\sigma\cat m}\sub
U_{\tau\cat m}\sub
U_{\tau\cat n}.
\qedhere
\]
\epf

For each point $x\in X$, we define a function $f_x\in\NN$ by induction on its argument $n$:
\[
f_x(n):=\min\set{m}{x\in U_{(f_x(1),\dotsc,f_x(n-1))}}.
\]

\bclm
The assignment $x\mapsto f_x$ is upper continuous.
\eclm
\bpf
Fix natural numbers $n$ and $m$.
We must show that the set $U:=\set{x}{f_x(n)\le m}$ is open.
Let $x$ be a point in $U$. It suffices to show that the set $U$ contains an open neighborhood of the point $x$.
Indeed, let
\[
V:=U_{(f_x(1))}\cap U_{(f_x(1),f_x(2))}\cap\dotsb\cap U_{(f_x(1),f_x(2),\dotsc,f_x(n))}.
\]
For each point $y\in V$, we have $y\in U_{(f_x(1))}$, and thus $f_y(1)\le f_x(1)$.
Since $y\in U_{(f_x(1),f_x(2))}$, we have $y\in U_{(f_y(1),f_x(2))}$ (Claim~\ref{clm:mon}),
and thus $f_y(2)\le f_x(2)$. Continuing this way, we have $f_y(n)\le f_x(n)\le m$.
\epf

Since the space $X$ is Menger, there is a function $g\in\NN$ with $f_x\lei g$ for each point $x\in X$.
Then Bob wins the play
\[
(
\cU_{()},U_{(g(1))},\cU_{(g(1))},U_{(g(1),g(2))},\dotsc
).
\]
Indeed, fix a point $x\in X$. Let $n$ be minimal with $f_x(n)\le g(n)$.
Then, by Claim~\ref{clm:mon}, we have
$x\in U_{( f_x(1),\dotsc,f_x(n-1),f_x(n) )}\sub U_{( g(1),\dotsc,g(n-1),g(n) )}$.
This completes the proof of Hurewicz's Theorem.
The proof establishes the following lemma.

\blem
\label{lem:sting}
Fix a strategy for Alice, consisting of countable increasing covers that get finer as the initial sequence increases. Define the corresponding functions $f_x$. If $n$ is minimal with
$f_x(n)\le g(n)$, then $x\in U_{( g(1),\dotsc,g(n-1),g(n) )}$.\qed
\elem

We prove the present theorem by applying Lemma~\ref{lem:sting} at all possible nodes
of the strategy tree.
For each natural number $k$ and each finite sequence $\sigma\in\bbN^k$,
define a sequence $f^\sigma_x\in\FinSeqs$ by
$(f^\sigma_x(1),\dotsc,f_x(k)):=\sigma$, and
$f^\sigma_x(n):=\min\set{m}{x\in U_{(f_x(1),\dotsc,f_x(n-1))}}$
for $n>k$.

By Claim~\ref{clm:mon}, the set
$\set{f^\sigma_x}{x\in X}$ is not dominating.
It follows that the countable union
\[
\bigcup_{\sigma\in\FinSeqs}\set{f^\sigma_x}{x\in X}=\set{f^\sigma_x}{x\in X,\sigma\in\FinSeqs}
\]
is not dominating. Let $g\in\NN$ be a witness for that.

Fix a point $x\in X$. Let $n_1$ be minimal with $f^{()}_x(n_1)\le g(n_1)$.
By Lemma~\ref{lem:sting}, we have
$x\in U_{( g(1),\dotsc,g(n_1) )}$.
Let $n_2>n_1$ be minimal with
$f^{( g(1),\dotsc,g(n_1) )}_x(n_2)\le g(n_2)$.
Applying Lemma~\ref{lem:sting} to the strategy beginning at the node
$( g(1),\dotsc,g(n_1) )$, we have $x\in U_{( g(1),\dotsc,g(n_2) )}$, and so on.
This completes the proof of the theorem.
\epf

\ed